\documentclass[sn-mathphys]{sn-jnl}

\jyear{2022}%

\usepackage{amssymb,latexsym, amsthm,epsf}
\newif\ifXY\XYtrue
\ifXY
\input xy
\usepackage{xy}

\xyoption{all} %
\xyoption{curve} %
\xyoption{arc} \fi 

\usepackage{xy}

\theoremstyle{thmstyletwo}%
\newtheorem{lemma}{Lemma}[section]%
\newtheorem{construction}{Construction}[section]%
\theoremstyle{thmstylethree}%
\newtheorem{definition}{Definition}[section]%
\newtheorem{corollary}{Corollary}[section]%
\newtheorem{theorem}{Theorem}[section]%
\usepackage{amscd,amsmath}
\usepackage{enumitem}

\def\C{{\mathbb C}}

\def\P{{\mathbb P}}

\def \ta{\tau}

\def \ta1{\tau_1}

\def \G{\Gamma}

\def \do{\dots}

\newcommand{\Xgal}{X_{Gal}}

\newcommand{\ug}[1]{\Gamma_#1}
\newcommand{\Ggal}{\widetilde{G}}

\makeatletter
\newcommand{\uGammaSq}[2]{\begin{equation}\label{#1}
	\ug{#2}^2\uGammaSqChecknextarg}
\newcommand{\uGammaSqChecknextarg}{\@ifnextchar\bgroup{\uGammaSqGobblenextarg}{ = e,
		\end{equation} }}
\newcommand{\uGammaSqGobblenextarg}[1]{ = \ug{#1}^2\@ifnextchar\bgroup{\uGammaSqGobblenextarg}{  = e
\end{equation}}}
\makeatother

\newcommand{\todo}[1]{\textcolor{red}{[TODO: #1]}}

\newcommand{\ubegineq}[1]{\begin{equation}\label{#1}}
\newcommand{\uendeq}{\end{equation}}
\newcommand{\trip}[2] {{\langle \ug{#1}, \ug{#2} \rangle}}

\newcommand{\comm}[2] {{[\ug{#1}, \ug{#2}]}}

\newcommand{\uThreePointInner}[5]{
\begin{equation}\label{#5-1}
\langle\Gamma_{#2'},\Gamma_{#3}\rangle=\langle\Gamma_{#2'},\Gamma_{#3'}\rangle=
\langle\Gamma_{#2'},\Gamma_{#3}^{-1}\Gamma_{#3'}\Gamma_{#3}\rangle=e
\end{equation}
\begin{equation}\label{#5-2}
\langle\Gamma_{#2},\Gamma_{#3}\rangle=\langle\Gamma_{#2},\Gamma_{#3'}\rangle=
\langle\Gamma_{#2},\Gamma_{#3}^{-1}\Gamma_{#3'}\Gamma_{#3}\rangle=e
\end{equation}
\begin{equation}\label{#5-3}
\Gamma_{#3'}\Gamma_{#3}\Gamma_{#2'}\Gamma_{#3}\Gamma_{#2'}^{-1}\Gamma_{#3}^{-1}\Gamma_{#3'}^{-1} = \Gamma_{#4'}
\end{equation}
\begin{equation}\label{#5-4}
\Gamma_{#3'}\Gamma_{#3}\Gamma_{#2'}\Gamma_{#3'}\Gamma_{#2'}^{-1}\Gamma_{#3}^{-1}\Gamma_{#3'}^{-1} = \Gamma_{#4'}\Gamma_{#4}\Gamma_{#4'}^{-1}
\end{equation}
\begin{equation}\label{#5-5}
\Gamma_{#3'}\Gamma_{#3}\Gamma_{#2}\Gamma_{#3}\Gamma_{#2}^{-1}\Gamma_{#3}^{-1}\Gamma_{#3'}^{-1} = \Gamma_{#4'}
\end{equation}
\begin{equation}\label{#5-6}
\Gamma_{#3'}\Gamma_{#3}\Gamma_{#2}\Gamma_{#3'}\Gamma_{#2}^{-1}\Gamma_{#3}^{-1}\Gamma_{#3'}^{-1} = \Gamma_{#4'}\Gamma_{#4}\Gamma_{#4'}^{-1}.
\end{equation}
}

\newcommand{\uFourPointInner}[6]{
	\begin{equation}\label{#6-1}
\langle\Gamma_{#2'},\Gamma_#3\rangle=\langle\Gamma_{#2'},\Gamma_{#3'}\rangle=\langle\Gamma_{#2'},\Gamma_#3^{-1}\Gamma_{#3'}\Gamma_#3\rangle=e
	\end{equation}
	\begin{equation}\label{#6-2}
	\langle\Gamma_#4,\Gamma_#5\rangle=\langle\Gamma_{#4'},\Gamma_#5\rangle=\langle\Gamma_#4^{-1}\Gamma_{#4'}\Gamma_#4,\Gamma_#5\rangle=e
	\end{equation}
	\begin{equation}\label{#6-3}
	[\Gamma_{#3'}\Gamma_#3\Gamma_{#2'}\Gamma_#3^{-1}{\Gamma_{#3'}}^{-1},\Gamma_#5] = e
	\end{equation}
	\begin{equation}\label{#6-4}
	 [\Gamma_{#3'}\Gamma_#3\Gamma_{#2'}\Gamma_#3^{-1}{\Gamma_{#3'}}^{-1},\Gamma_#4^{-1}{\Gamma_{#4'}}^{-1}\Gamma_#5^{-1}\Gamma_{#5'}\Gamma_#5\Gamma_{#4'}\Gamma_#4] = e
	\end{equation}
	\begin{equation}\label{#6-5}
	\langle\Gamma_#2,\Gamma_#3\rangle=\langle\Gamma_#2,\Gamma_{#3'}\rangle=\langle\Gamma_#2,\Gamma_#3^{-1}\Gamma_{#3'}\Gamma_#3\rangle=e
	\end{equation}
	\begin{equation}\label{#6-6}
\langle\Gamma_#4,\Gamma_#5^{-1}\Gamma_{#5'}\Gamma_#5\rangle=\langle\Gamma_{#4'},\Gamma_#5^{-1}\Gamma_{#5'}\Gamma_#5\rangle=\langle\Gamma_#4^{-1}
\Gamma_{#4'}\Gamma_#4,\Gamma_#5^{-1}\Gamma_{#5'}\Gamma_#5\rangle=e
	\end{equation}
	\begin{equation}\label{#6-7}
	[\Gamma_{#3'}\Gamma_#3\Gamma_#2\Gamma_#3^{-1}{\Gamma_{#3'}}^{-1},\Gamma_#5^{-1}\Gamma_{#5'}\Gamma_#5] = e
	\end{equation}
	\begin{equation}\label{#6-8}
	[\Gamma_{#3'}\Gamma_#3\Gamma_#2\Gamma_#3^{-1}{\Gamma_{#3'}}^{-1}, \Gamma_#4^{-1}{\Gamma_{#4'}}^{-1}\Gamma_#5^{-1}{\Gamma_{#5'}}^{-1}\Gamma_#5\Gamma_{#5'}\Gamma_#5\Gamma_{#4'}\Gamma_#4] = e
	\end{equation}
	\begin{equation}\label{#6-9}
	\Gamma_{#3'}\Gamma_#3\Gamma_{#2'}\Gamma_#3\Gamma_{#2'}^{-1}\Gamma_#3^{-1}{\Gamma_{#3'}}^{-1} = \Gamma_#5\Gamma_{#4'}\Gamma_#5^{-1}
	\end{equation}
	\begin{equation}\label{#6-10}
	\Gamma_{#3'}\Gamma_#3\Gamma_{#2'}\Gamma_{#3'}\Gamma_{#2'}^{-1}\Gamma_#3^{-1}{\Gamma_{#3'}}^{-1} = \Gamma_#5\Gamma_{#4'}\Gamma_#4{\Gamma_{#4'}}^{-1}\Gamma_#5^{-1}
	\end{equation}
	\begin{equation}\label{#6-11}
	\Gamma_{#3'}\Gamma_#3\Gamma_#2\Gamma_#3\Gamma_#2^{-1}\Gamma_#3^{-1}{\Gamma_{#3'}}^{-1} = \Gamma_#5^{-1}\Gamma_{#5'}\Gamma_#5\Gamma_{#4'}\Gamma_#5^{-1}{\Gamma_{#5'}}^{-1}\Gamma_#5
	\end{equation}
	\begin{equation}\label{#6-12}
	\Gamma_{#3'}\Gamma_#3\Gamma_#2\Gamma_{#3'}\Gamma_#2^{-1}\Gamma_#3^{-1}{\Gamma_{#3'}}^{-1} = \Gamma_#5^{-1}\Gamma_{#5'}\Gamma_#5\Gamma_{#4'}\Gamma_#4\Gamma_{#4'}^{-1}\Gamma_#5^{-1}{\Gamma_{#5'}}^{-1}\Gamma_#5.
	\end{equation}
}

\newcommand{\uFivePointOuter}[7]{
	\todo{}
}

\newcommand{\uParasit}[3]{
\begin{equation}\label{#3}
[\Gamma_{#1},\Gamma_{#2}] = [\Gamma_{#1'},\Gamma_{#2}]  = [\Gamma_{#1},\Gamma_{#2'}]  = [\Gamma_{#1'},\Gamma_{#2'}] = e
\end{equation}}

\makeatletter
\newcommand{\uProjRel}[2]{
\begin{equation}\label{#1}
\Gamma_{#2'}\Gamma_{#2}\uProjRelChecknextarg}
\newcommand{\uProjRelChecknextarg}{\@ifnextchar\bgroup{\uProjRelGobblenextarg}{ = e
	\end{equation}}}
\newcommand{\uProjRelGobblenextarg}[1]{\Gamma_{#1'}\Gamma_{#1}\@ifnextchar\bgroup{\uProjRelGobblenextarg}{  = e. \end{equation}}}
\makeatother

\begin{document}

\title[Non-planar degenerations and related fundamental groups] {Non-planar degenerations and related fundamental groups}


\author[]{{Meirav Amram}} 

\affil[]{\orgdiv{Department of Mathematics}, \orgname{SCE}, \orgaddress{\street{Jabotinsky 84}, \city{Ashdod}, \country{Israel}}
\\ {ORCID: 0000-0003-4912-4672 \ \ \ Email: meiravt@sce.ac.il}}


\abstract
{We present a preliminary investigation of algebraic surfaces that have non-planar degenerations, along with their Galois covers and fundamental groups.

 Specifically, we investigate the tetrahedron and the double tetrahedron. 
 The resulting fundamental groups indicate that the tetrahedron and the double tetrahedron are in different components of the moduli space of algebraic surfaces.

 Non-planar degenerations provide a way to understand the behavior of algebraic surfaces in a broader context. By considering degenerations that deviate from the idealized planar case, we gain insights into the structures and properties of surfaces in more general settings. This helps develop a more comprehensive understanding of algebraic surfaces and their behavior.
 
This study aims to determine the fundamental groups of the Galois covers of some algebraic surfaces; these groups are invariants of the classification of surfaces in the moduli space of algebraic surfaces. 
 Our findings can advance the classification of surfaces and provide further links between algebraic geometry, group theory,  and the topology of degenerative processes and their properties.}




\maketitle

\section{Introduction}\label{outline}
Classifying algebraic surfaces and studying their moduli space are some of the most thoroughly investigated subjects in algebraic geometry and topology (see e.g., Catanese \cite{C1,C2}). Algebraic surfaces, as geometric objects, can be investigated via their convexity and curvature \cite{Betti}, their degenerations  \cite{CCFM,2008}, their related moduli properties \cite{Lid}, Chern classes \cite{Ma}, fundamental groups \cite{Alex1,Li08}, and  many other geometrical and topological invariants.

One of the known invariants of classification is the fundamental group $\pi_1(\Xgal)$ of the Galois cover $\Xgal$ of an algebraic surface $X$, with respect to a generic projection $f$ to the projective plane $\C\P^2$. The group $\pi_1(\Xgal)$ has geometric significance in the classification of surfaces because it is isomorphic for all the surfaces in the same connected components of the moduli space of algebraic surfaces.

Non-planar degenerations appear in the study of moduli spaces of algebraic surfaces, where the moduli space of all algebraic surfaces of a given type is parameterized.  The planar degenerations are those that can be depicted on a piece of paper; they are defined in  Definition \ref{def_planar_representation}. Non-planar degenerations arise naturally when we identify edges in planar degenerations and then glue them together, which induces appearances of high multiple singularities, and therefore makes computations much more difficult than before. 
Non-planar degenerations have important applications in various fields, such as algebraic geometry, topology, and physics.  By studying non-planar degenerations, we can better understand the geometry and topology of these moduli spaces of algebraic surfaces, leading to important advancements in these fields.

The motivation to study non-planar degenerations of algebraic surfaces lies in gaining a deeper understanding of surface geometry and invariants, as well as their applications in related fields. 
 We hope to combine our future results for $\pi_1(\Xgal)$ with the results obtained in \cite{AG} and look for their generalizations.

Examples of surfaces that have non-planar degenerations appear in \cite{pillow,AG,CP1*T,T*T}; these are products of a complex torus with the projective line, a product of a complex torus with a complex torus, and $K_3$ surfaces that are known to have {\it pillow degenerations}.  In \cite{AG}, we used the  Reidemeister-Schreier method to find  $\pi_1(\Xgal)$ where $X=\mathbb{CP}^1 \times T$  (and $T$ is a complex torus). In \cite{pillow}, we presented only the list of braids for a certain $K3$ surface, but we did not find fundamental groups.
Other previous works deal primarily with planar degenerations, for example, in \cite{6degree}, we give a full classification of 29 cases of degree 6 planar degenerations.

In \cite{8degree} we carried out the study of a non-planar degeneration of degree 8, which is isomorphic to an octahedron and has singularities of multiplicity 4.  Our goal is to obtain general results about fundamental groups related to  non-planar degenerations, and this study will be much more challenging. We suggest a new direction of research that will address even non-planar degenerations, from the simplest to the most complex, including determining the fundamental group in the extensive algebraic use of group theory and computational methods.

In this paper, we determine  $\pi_1(\Xgal)$ for the objects $T_{(4)}$ and $D(T_{(4)})$ (a tetrahedron and a double tetrahedron, respectively), as described in Theorems \ref{TetGal} and \ref{DoubleTetGal}. These surfaces have non-planar degenerations of the smallest degrees 4 and 6, respectively.  As their $\pi_1(\Xgal)$ are not isomorphic, then these two surfaces are in different connected components of the moduli space.

To find $\pi_1(\Xgal)$, we use a generic projection of a surface $X$ and its degeneration $X_0$ onto $\mathbb{CP}^2$ to get their branch curves $S$ and $S_0$, respectively. We can then use the van Kampen Theorem \cite{vk} to get a presentation of the fundamental group $\pi_1(\mathbb{CP}^2-S)$ of the complement of $S$. We can also construct the dual graph that corresponds to the group and to its presentation. 
The primary advantage of this algorithm is that we use it when the branch curve $S$ and the fundamental group $\pi_1(\mathbb{CP}^2-S)$ are difficult to describe, especially because singularities with high multiplicities appear. This algorithm recovers the curve $S$ and a significant amount of geometric and algebraic information, giving us the group $\pi_1(\Xgal)$.
 In \cite{dettweiler} there is an interesting explanation about the correspondence of plane curves and groups. The works \cite{Deg1,Deg2,Iten} connect the topics of degenerations, topology of plane curves, and fundamental groups, and they are also a basis for understanding singular points in curves.

This paper is organized as follows:  In Section \ref{section:method}, we present the algorithm and give notations to the fundamental groups. We provide notations for the generators and give the presentations of the groups via the van Kampen Theorem. In Section \ref{sec3}, we explain the degenerations of interest and construct them.
In Subsections \ref{Tetra-sec} and \ref{DoubleTetra-sec}, we calculate the fundamental groups of the Galois covers related to the tetrahedron $T_{(4)}$ and the double tetrahedron $D(T_{(4)})$, respectively. In addition, in Subsection \ref{DoubleTetra-sec} we give some necessary details on Coxeter groups and dual graphs. In Section \ref{conc}, we present the summary of the results and the conclusion.

\section{Preliminaries and details of the algorithm}\label{section:method}
In this section, we provide some basic background and outline process steps. We recommend  \cite{2008} as a source for relevant background about degenerations of projective structures. Calabri-Ciliberto-Flamini-Miranda and Ciliberto-Lopez-Miranda are among the pioneers in research on degenerations and from their works \cite{CCFM,CLM} one can get a basic understanding of this subject. 
The book \cite{GoWed} gives us an excellent background about projective varieties and affine algebraic sets along with their topological properties. The book  \cite{Barth} is a cornerstone for understanding algebraic surfaces and their properties. Brazas \cite{brazas},  Auroux-Donaldson-Katzarkov-Yotov \cite{ADKY}, and Amram-Friedman-Teicher \cite{AFT09,AFT03} supply some background and examples on fundamental groups.

We explain the algorithm of degeneration and regeneration that provides a fundamental group of the Galois cover of an algebraic surface embedded in projective space $\mathbb{CP}^n $.
	We take a generic projection of the surface $X$ onto the projective plane $\mathbb{CP}^2$.
	We get the branch curve $S$ in $\mathbb{CP}^2$. The branch curve $S$ indicates singularities and the overall shape of the surface $X$.  Since $S$ is very difficult to describe, we will use the process of degeneration, as indicated in Definition \ref{deg-def}.
	
	\begin{definition}\label{deg-def}
		A degeneration of $X$ is a proper surjective morphism with connected fibers
$\pi  : V \rightarrow \mathbb{C}$ from an algebraic variety $V$, such that the restriction $\pi:
V \setminus \pi^{-1}(0) \rightarrow \mathbb{C} \setminus \{0\}$
is smooth, and that $\pi^{-1}(t) \cong X$ for $t\neq0$.

When $X$ is projective with an embedding $\rho: X \hookrightarrow \mathbb{CP}^n$,
 a degeneration of $X$ is called a projective degeneration
of $\rho$ if there exists a morphism $F  : \ V \rightarrow
\mathbb{CP}^n \times \mathbb{C}$ such that the restriction
$\ F_t = F \mid_{\pi^{-1}(t)} \ : \pi^{-1}(t) \rightarrow \mathbb{CP}^n
\times \{t\}$ is an embedding of $\pi^{-1}(t)$ for all $t\in \mathbb{C}$
and that $F_1 = \rho$ under the identification of $\pi^{-1}(1)$ with $X$.
	\end{definition}
	
The first step of the algorithm is to construct a flat degeneration $X_0$ of $ X $, as a union of planes. Each plane in $X_0$ is linearly, algebraically, and analytically isomorphic to the projective plane $\mathbb{CP}^2$.     
We consider degenerations with only two planes intersecting at an edge (because this is the generic case). 
The branch curve $S_0$ of $X_0$ is a union of lines that  are the projections of those edges; it has singularities of the type defined below.

\begin{figure}[H]
\begin{center}
\scalebox{0.7}{\includegraphics{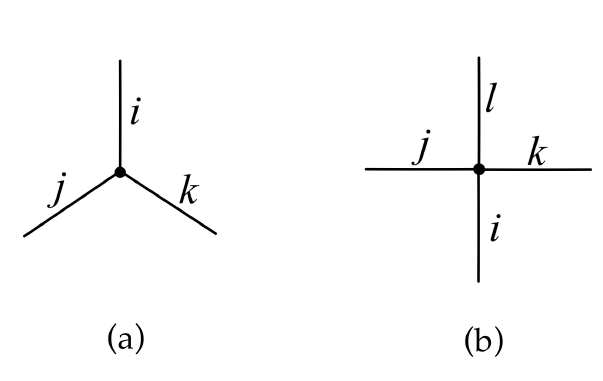}}
\end{center}
\setlength{\abovecaptionskip}{-0.15cm}
\caption{Singularities of multiplicities $3$ and $4$}\label{3-4point}
\end{figure}

\begin{definition}	\label{inner-pt}
	An inner singularity of multiplicity $k$ is an intersection point of $k$ planes, such that two neighboring planes share a direct edge at their intersection, while non-neighboring planes touch each other only at the singularity.
	\end{definition}
In this paper we have $k=3,4$, see Fig. \ref{3-4point}.
Note that similar singularities appear in \cite{6degree,CCFM}, and it can be helpful to use the information there as well.

One of the main tools that is embedded in the algorithm is a reverse process of degeneration, which is called a {regeneration}.
	We recover $ S $ from $ S_0 $ by regeneration lemmas from \cite{BGT2}. The degeneration and regeneration process is illustrated in the following diagram.

	\[\begin{CD}
	X\subseteq \mathbb{CP}^n  @>\text{degeneration}>> X_0\subseteq \mathbb{CP}^n \\
	@V\text{generic~ projection}VV                      @VV\text{generic~ projection}V \\
	S\subseteq \mathbb{CP}^2 @<\phantom{regeneration}<\text{regeneration}< S_0\subseteq \mathbb{CP}^2
	\end{CD}\]
	
	\vspace{0.2cm}	
	
	The regeneration process is the second step of the algorithm; a line in $S_0$ regenerates to a conic, and  the tangency point  we get during this process will regenerate later to three cusps, see \cite{19} for more details.
Therefore, the regenerated branch curve $S$ is cuspidal, and its degree is double that of~$S_0$.


Let $G:=\pi_1(\mathbb{CP}^2-S)$ be the fundamental group of the complement of curve $S$. The standard generators in the van Kampen presentation of this group are $\G_1, \G_{1'}, \dots, \G_{m},\G_{m'}$, where $m$ denotes the number of lines in the degenerated curve $S_0$ (or equivalently, the number of edges in $X_0$).

To get a presentation of $G$ using the generators $\{\G_{j}, \G_{j'}\}$ ($1 \leq j \leq m$), we need the van Kampen Theorem \cite{vk} with (a)-(d), as shown in Theorem \ref{vkthm}:
\begin{theorem}\label{vkthm}
	\begin{enumerate}
		\item [(a)] For a branch point of a conic, we have the relation $\G_{j} = \G_{j'}$.
		\item [(b)] For nodes, we have  $[\G_{i},\G_{j}]=\G_{i}\G_{j}\G_{i}^{-1}\G_{j}^{-1}=e$.
		\item [(c)] For cusps,  we have  $\langle\G_{i},\G_{j}\rangle=\G_{i}\G_{j}\G_{i}\G_{j}^{-1}\G_{i}^{-1}\G_{j}^{-1}=e$.
\item [(d)] To the presentation we add the projective relation $\prod\limits_{j=m}^1 \G_{j'}\G_{j}=e$.
\item [(e)]  In addition to the van Kampen presentation, we add commutations, which come from lines in degenerations that do not ordinarily intersect, but when projecting each one of them onto $\C\P^2$, they will (further details about these commutations are provided in \cite{MoTe87}). 
	\end{enumerate}
\end{theorem}

We recall from \cite{MoTe87} that if
	$f : X\rightarrow \mathbb{CP}^2$ is a generic projection of degree $n$,
then $\Xgal$, the Galois cover, is defined as follows:
	$$\Xgal=\overline{(X \times_{\mathbb{CP}^{2}}\dots \times_{\mathbb{CP}^{2}} X)-\triangle},$$
	where the product is taken $n$ times and the diagonal 
 $\triangle$ is the set of elements $(x_1,\dots,x_k)\in X^k$ such that $x_i=x_j$ for some $i\ne j$. 
We define $$ \Ggal:=\frac{G}{\langle \Gamma_j^2, \Gamma_{j'}^2\rangle},$$ and use the exact sequence from \cite{MoTe87}:
	\begin{equation}\label{M-T}
	0 \rightarrow \pi_1(\Xgal) \rightarrow \widetilde{G} \rightarrow S_n \rightarrow 0,
	\end{equation}
	where the second map takes the generators $ \Gamma_j $ and $ \Gamma_{j'} $ of $G$ to the transposition of the two planes that intersect at line $ j $.
	We thus obtain a presentation of the fundamental group $\pi_1(\Xgal)$ of the Galois cover.

Galois covers of algebraic surfaces were investigated by  Liedtke \cite{Li08}, Moishezon-Teicher \cite{MoTe87}, and Gieseker \cite{Gie}.
Fundamental groups of Galois covers were studied for $\mathbb{CP}^1 \times T$ (where $T$ is a complex torus) in Amram-Goldberg \cite{AG} and Amram-Tan-Xu-Yoshpe \cite{CP1*T}, and for toric varieties in Amram-Ogata \cite{Ogata}. Furthermore, the groups were studied for surfaces with Zappatic singularity of type $E_k$ in Amram-Gong-Tan-Teicher-Xu \cite{ZAPP}, and for surfaces with degenerations of degree $6$ in Amram-Gong-Sinichkin-Tan-Xu-Yoshpe \cite{6degree}. Grigorchuk \cite{Gri} provides a great deal of information about these groups and their associated Coxeter groups. 
 Moreover, in \cite{gan}, there are explanations about the correspondence between graphs and fundamental groups, and in \cite{Gon}, the authors work with non-isomorphic fundamental groups and Galois groups.

\section{The degenerations of $T_{(4)}$ and $D(T_{(4)})$}\label{sec3}

In this paper we consider the degenerations of the tetrahedron and the double tetrahedron. Before we define the class of degenerations to which those examples correspond, we recommend  \cite{saidi} as a source for relevant information about some combinatorial degeneration data and related Galois covers. 

\begin{figure}[ht]
\begin{center}
\scalebox{0.30}{\includegraphics{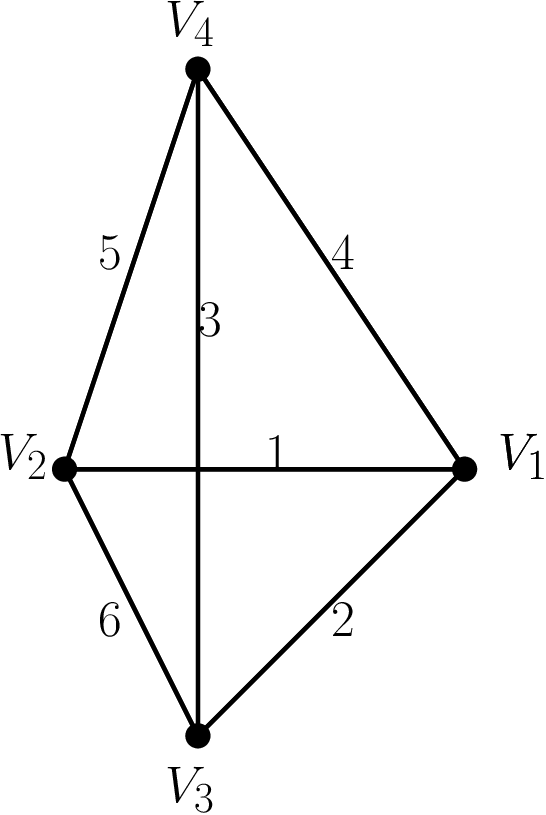}}
\end{center}
\setlength{\abovecaptionskip}{-0.15cm}
\caption{The degeneration of the tetrahedron $T_{(4)}$}\label{Tetrahedron}
\end{figure}

First we recall the definition of a planar degeneration, as shown in Definition \ref{def_planar_representation}:
\begin{definition}\label{def_planar_representation}
	A degeneration of smooth toric surface $ X $ into a union of planes $ X_0 $, is said to have a planar representation if:
	\begin{enumerate}
		\item
		No three planes in $ X_0 $ intersect in a line.
		\item
		There exists a simplicial complex with a connected interior embedded in $ \mathbb{R}^2 $, such that its triangles correspond bijectively to irreducible components of $ X_0 $, and the dimension of the intersection of any set of triangles is equal to the dimension (over $\mathbb{C}$) of the intersection of the corresponding components.
	\end{enumerate}
\end{definition}
While neither $ T_{(4)} $ nor $ D(T_{(4)}) $ are of this type, property (1) in Definition \ref{def_planar_representation} still holds for these degenerations.
We thus can define a more general set of examples, as shown in Definition \ref{def_combinatorially_homeomorphic}:
	\begin{definition}\label{def_combinatorially_homeomorphic}
		Let $ \Omega $ be a connected topological manifold of dimension $ m $.
		A degeneration of a smooth $m$-dimensional toric variety $ X $ into a union of projective spaces $ X_0 $ is said to be combinatorially homeomorphic to $ \Omega $ if:
		\begin{enumerate}
			\item
			No three irreducible components of $ X_0 $ intersect in a co-dimension 1 set.
			\item
			There exists a simplicial complex of pure dimension $ m $ with a connected interior that is homeomorphic to $ \Omega $ such that  its $m$-simplices correspond bijectively to irreducible components of $ X_0 $, and the dimension of the intersection of any set of triangles is equal to the dimension (over $\mathbb{C}$) of the intersection of the corresponding components.
		\end{enumerate}
	\end{definition}
This definition generalizes Definition \ref{def_planar_representation} because planar degenerations are precisely those that are combinatorially homeomorphic to a disc. Now, the degeneration of the tetrahedron $T_{(4)}$ is shown in Fig. \ref{Tetrahedron}.

\begin{figure}[H]
\begin{center}
\scalebox{0.60}{\includegraphics{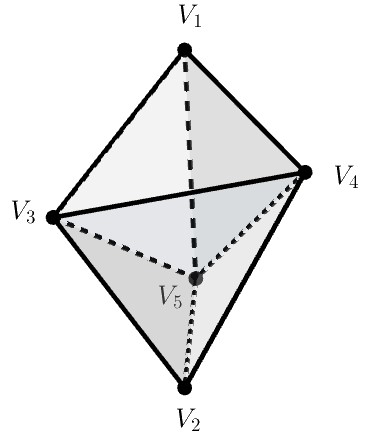}}
\end{center}
\setlength{\abovecaptionskip}{-0.15cm}
\caption{The degeneration of the double tetrahedron $D(T_{(4)})$}\label{DoubleTetra}
\end{figure}

The degeneration  of the double tetrahedron $D(T_{(4)})$ is shown in Fig. \ref{DoubleTetra}. 

\begin{figure}[ht]
\begin{center}
\scalebox{0.25}{\includegraphics{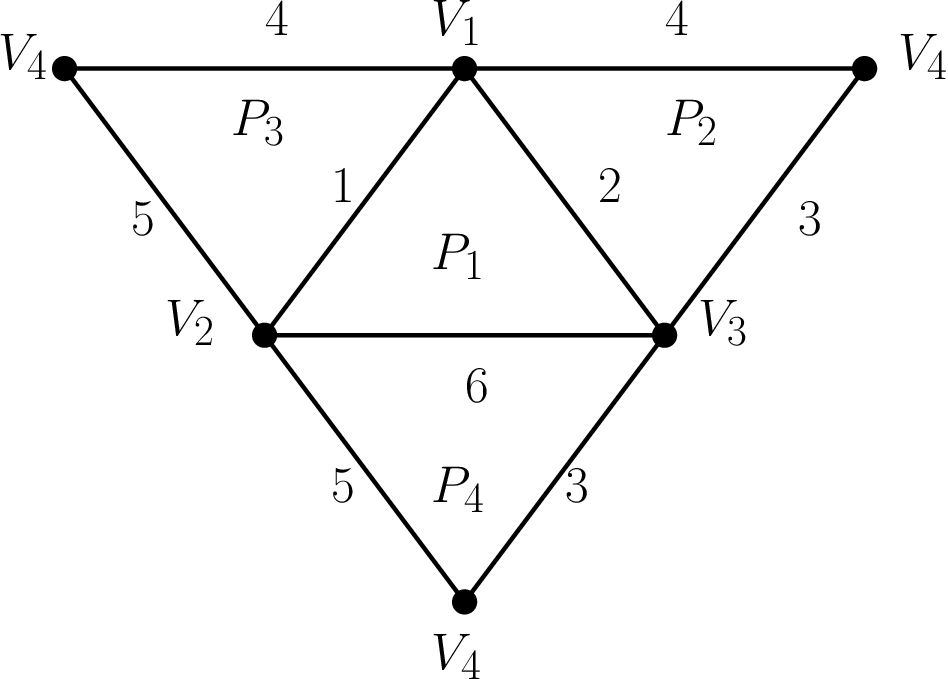}}
\end{center}
\setlength{\abovecaptionskip}{-0.15cm}
\caption{The tetrahedron $T_{(4)}$ before identification}\label{Tetrahedron-planarform}
\end{figure}

\begin{figure}[H]
\begin{center}
\scalebox{0.30}{\includegraphics{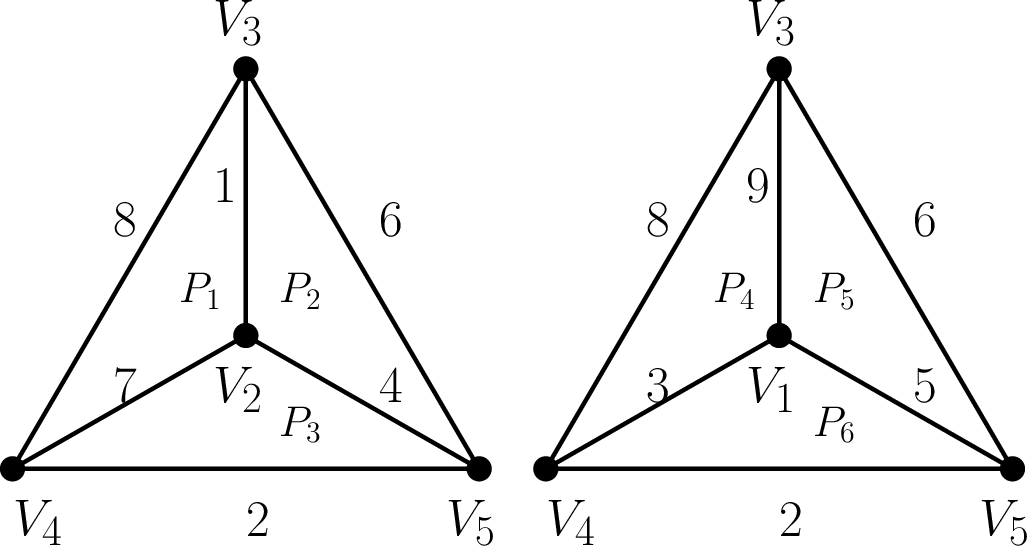}}
\end{center}
\setlength{\abovecaptionskip}{-0.15cm}
\caption{The double tetrahedron $D(T_{(4)})$ before identification}\label{twopieces}
\end{figure}

Both degenerations are non-planar, and we get them by identifications of common edges. For example, in Fig. \ref{Tetrahedron-planarform} we identify the two edges numbered by 3, the two edges numbered by 4, and the two edges numbered by 5. These identifications yield $T_{(4)}$ in Fig. \ref{Tetrahedron}. 
In a similar way, the glue of the two pieces along the edges having the same number in Fig. \ref{twopieces} yields   $D(T_{(4)})$ in Fig. \ref{DoubleTetra}. 
This is the reason that both $ T_{(4)} $ and $ D(T_{(4)}) $ are combinatorially homeomorphic to a sphere. 
In addition, we see that the pillow degeneration considered in \cite{pillow} is combinatorially homeomorphic to a sphere. In contrast, the degenerations of $ \mathbb{CP}^1\times T $ in \cite{AG} and of $ T\times T $ in \cite{TT} are combinatorially homeomorphic to a cylinder and a topological torus, respectively.

The lowest degree possible for non-planar degeneration $X_0$ of this kind is $4$,  for which we have the tetrahedron $T_{(4)}$.
This tetrahedron is a union of four planes, see Fig. \ref{Tetrahedron}. The smallest degree of non-planar degeneration of an algebraic surface is 4, as gluing two planar pieces  creates singularities that cannot be flattened into a plane.
The next in line with a higher degree is $D(T_{(4)})$; it is of degree 6, as shown in Fig. \ref{DoubleTetra}.  
 The double tetrahedron can be considered as a glue of two tetrahedrons without bases (they are glued along the edges of the missing bases). More details about both surfaces are discussed in the construction below. 
	Because we study those degenerations with the glue of two pieces that look the same, we get the even number of triangles.


\begin{construction}[$ T_{(4)} $]
We construct the simplicial complex corresponding to the tetrahedron $T_{(4)}$ that appears  in Fig. \ref{Tetrahedron-planarform}.
\begin{itemize}
  \item Four triangles correspond to four planes, denoted by $P_1, P_2, P_3, P_4$.
  \item Edges 1, \dots, 6 are the lines that are the intersections between the planes.
  \item Planes $P_1$ and $P_2$ have a common edge 2, $P_1$ and $P_3$ have a common edge 1, $P_1$ and $P_4$ have a common edge 6, $P_2$ and $P_3$ have a common edge 4, $P_2$ and $P_4$ have a common edge 3, and $P_3$ and $P_4$ have a common edge 5.
  \item $V_1, V_2, V_3, V_4$ are four vertices. Planes $P_1, P_2$, and $P_3$ meet in vertex $V_1$;  planes $P_1, P_3$, and $P_4$ meet in vertex $V_2$; planes $P_1, P_2$, and $P_4$ meet in vertex $V_3$; and planes $P_2, P_3$, and $P_4$ meet in vertex $V_4$.
\end{itemize}
We denote   $X_0 = \bigcup_{i=1}^4 P_i$.

\end{construction}

\begin{construction}[$ D(T_{(4)}) $]
The simplicial complex corresponding to the degeneration of $ D(T_{(4)}) $  can be seen in Fig. \ref{twopieces}.
\begin{itemize}
\item Six triangles correspond to six planes, denoted by $P_1, \dots, P_6$.
\item Edges 1, \dots, 9 are the lines that are the intersections between the planes.
\item Plane $P_1$ shares edges 1, 7, and 8 with planes $P_2$, $P_3$, and $P_4$, respectively.
Plane $P_2$ shares edges 4 and 6  with planes $P_3$ and $P_5$, respectively.
Plane $P_6$ shares edges 2, 3, and 5, with planes $P_3$, $P_4$, and $P_5$, respectively.
Planes $P_4$ and $P_5$ share edge 9.
\item $V_1, \dots, V_5$ are five vertices. Planes $P_4, P_5$, and $P_6$ meet in vertex $V_1$;  planes $P_1, P_2$, and $P_3$ meet in vertex $V_2$; planes $P_1, P_2, P_4$, and $P_5$ meet in vertex $V_3$; planes $P_1, P_3, P_4$, and $P_6$ meet in vertex $V_4$, and planes $P_2, P_3, P_5$, and $P_6$ meet in vertex $V_5$.
\end{itemize}
We denote   $X_0 = \bigcup_{i=1}^6 P_i$.
\end{construction}

To compute the fundamental groups of the Galois covers of these two surfaces,
we work with dual graphs, which were presented in \cite{Law} and \cite{RTV}. Each graph $ T $ relates to each $ X_0 $, and  is defined as follows:
The vertices of $ T $ are in bijection with the planes in $ X_0 $,
and the vertices corresponding to the planes $ P_i $ and $ P_j $ are connected by an edge
if $ P_i $ and $ P_j $ intersect in an edge. Each generator (or a pair of generators) in $G$ or in $\widetilde{G}$ represents an edge between two planes in the degeneration.

In $ T_{(4)} $ we have inner singularities of multiplicity $3$ (vertices $V_i$,  $i=1,\dots,4$ in Fig. \ref{Tetrahedron-planarform}). 
In $ D(T_{(4)}) $ we have inner singularities of multiplicity $3$ as well (vertices $V_1$ and $V_2$ in Fig. \ref{twopieces}). In addition, $ D(T_{(4)}) $ has singularities of multiplicity $4$ (vertices $V_i$, $i=3,4,5$).
Fig. \ref{3-4point}(a) depicts three lines, being globally ordered as $i < j < k$; they produce a multiplicity $3$ singularity. Fig. \ref{3-4point}(b) depicts an intersection of four lines, being globally ordered as $ i < j < k < l$.

 Our numerations of the edges in each of the degenerations  were  chosen so that locally we can write the relations in $G$ that come exactly from the singularities in Fig. \ref{3-4point}. Singularities with multiplicity $3$ were studied in detail in previous papers, such as  \cite{5degree} and \cite{6degree}; singularities of multiplicity $4$ were studied in \cite{6degree}, and  recently in \cite{8degree}, as well.  
 
 Lemmas \ref{lemma-3} (from \cite[Lemma 3.3]{6degree}) and \ref{lemma-4} (from \cite[Theorems 21,24]{Ogata}) summarize the relations in $G$ that are contributed by those singularities. We recall from Theorem \ref{vkthm} that relation $\langle\G_{i},\G_{j}\rangle=e$ means $\G_{i}\G_{j}\G_{i}=\G_{j}\G_{i}\G_{j}$.

\begin{lemma}\label{lemma-3}
     A point that is an intersection of three lines, as in Fig. \ref{3-4point}(a), contributes to $G$ the following list of relations:
     \uThreePointInner{p}{i}{j}{k}{3pt-rels}
\end{lemma}

\begin{lemma}\label{lemma-4}
 A point that is an intersection of four lines, as in Fig. \ref{3-4point}(b), contributes to $G$ the following list of relations:
\uFourPointInner{q}{i}{j}{k}{l}{4pt-rels}
\end{lemma}

In addition to the above relations that we get from the singularities, the following lemma corresponds to Lemma \ref{lemma-3} and speaks about a specific condition for a multiplicity $3$ singularity. If that specific condition holds, this lemma (which was proven in \cite{6degree}) can assist us in the computations of $\widetilde{G}$.

\begin{lemma}\label{3pt-in-bigmid}
	Given a singularity of multiplicity $3$ in $ X_0 $ with lines $ i < j < k $ as in Fig. \ref{3-4point}(a),  
if either $ \Gamma_j=\Gamma_{j'} $ or $ \Gamma_k=\Gamma_{k'} $ holds in $ \widetilde{G} $
then $ \Gamma_l=\Gamma_{l'} $ for all $ l\in \{i,j,k\} $.
	Moreover, for each $i$, which fulfills the condition  $ i < j < k $ in such singularity, the relation $ \Gamma_i=\Gamma_{i'} $ always holds. 

\end{lemma}




\bigskip

\section{Computing the fundamental group $\pi_1(\Xgal)$}\label{results}
In this section, we determine $\pi_1(\Xgal)$ for both surfaces $T_{(4)}$ and $D(T_{(4)})$. We apply the algorithm and methods described in detail in Section \ref{section:method}.

\subsection{The tetrahedron $T_{(4)}$}\label{Tetra-sec}

In this subsection, we consider the tetrahedron $T_{(4)}$ from Fig. \ref{Tetrahedron}.
The tetrahedron acquires a non-planar degeneration after the identifications along edges 3, 4, and 5 in Fig. \ref{Tetrahedron-planarform}.

In the proof of Theorem \ref{TetGal}, we calculate the fundamental group $\pi_1(\Xgal)$ 
of the Galois cover of the tetrahedron $T_{(4)}$.
\begin{theorem}\label{TetGal}
The fundamental group $\pi_1(\Xgal)$ of the Galois cover of the tetrahedron $ T_{(4)} $ is trivial.
\end{theorem}

\begin{verbatim}
Proof.
\end{verbatim}
There are six edges in $ T_{(4)} $ in Fig. \ref{Tetrahedron}. Projecting $ T_{(4)} $ onto $\mathbb{CP}^2$, we get curve $S_0$, which is an arrangement of six lines that are the projections of the six edges. Vertices $V_i$, $i=1,\dots,4$, are singularities of multiplicity~$3$. 
In the regeneration the degree is doubled, therefore the branch curve $S$ has degree 12. The generators in group $G$ are $\Gamma_{1},\Gamma_{1'},\dots,\Gamma_6,\Gamma_{6'}$.  Because we are interested in fundamental groups, we provide first the relations in~$G$ that come from each vertex. We use Lemma \ref{lemma-3} to write the relations induced by vertex $V_1$. 
\uThreePointInner{1}{1}{2}{4}{Tetra,vert1}
Similarly, we have relations induced by vertices $V_2$ (with indices $i=1,j=5,k=6$), $V_3$ (with indices $i=2,j=3,k=6$), and $V_4$ (with indices  $i=3,j=4,k=5$).

Then we collect all the above relations, together with the commutative relations of type ({e}) from Theorem \ref{vkthm} (see them in \eqref{Tetra,parasit-1-3}-\eqref{Tetra,parasit-4-6}) and the projective relation (in \eqref{Tetra,proj}).
\uParasit{1}{3}{Tetra,parasit-1-3}
\uParasit{2}{5}{Tetra,parasit-2-5}
\uParasit{4}{6}{Tetra,parasit-4-6}
\uProjRel{Tetra,proj}{6}{5}{4}{3}{2}{1}

Because we have the exact sequence \eqref{M-T} and we are interested in the group $\pi_1(\Xgal)$, we calculate now $\Ggal$. 
An immediate conclusion from Lemma \ref{3pt-in-bigmid} gives us  
$ \Gamma_1=\Gamma_{1'}$ in vertices $V_1$ and $V_2$, $ \Gamma_2=\Gamma_{2'}$ in vertex $V_3$,
and $ \Gamma_3=\Gamma_{3'}$ in vertex $V_4$. 
Then, we apply Lemma \ref{3pt-in-bigmid} again, this time with these equalities, and we get   
that $\Gamma_i=\Gamma_{i'}$ for $i= 1, \dots, 6$.
 These equalities simplify the relations of type (c) from Theorem \ref{vkthm} to
\begin{align}	
	\langle \Gamma_{1},\Gamma_{2}\rangle = \langle \Gamma_{1},\Gamma_{5}\rangle =
\langle \Gamma_{2},\Gamma_{3}\rangle = \langle \Gamma_{3},\Gamma_{4}\rangle =e, \label{trip}
	\end{align}	
the branch relations of type (a) from Theorem \ref{vkthm} to 
\begin{align}	
 \ug{4}=\ug{1}\ug{2}\ug{1}, \ \ \ug{5}=\ug{3}\ug{4}\ug{3}, \ \ \ug{6}=\ug{2}\ug{3}\ug{2}, \ \  \ug{6}=\ug{1}\ug{5}\ug{1},\label{equal}
	\end{align}	
and the commutations of type (b) from Theorem \ref{vkthm} to
$$\comm{1}{3} = \comm{2}{5} = \comm{4}{6} = e.$$

We continue the simplification of $\Ggal$ as follows: we substitute equalities from (\ref{equal}) inside relations from (\ref{trip}) to get more relations of type (c). For example, we substitute $\ug{5}=\ug{1}\ug{6}\ug{1}$ in $\trip{1}{5}=e$ and get $\trip{1}{6}=e$; we substitute $\ug{2}=\ug{1}\ug{4}\ug{1}$ in $\trip{1}{2}=e$ and get $\trip{1}{4}=e$; then we substitute  $\ug{1}=\ug{2}\ug{4}\ug{2}$ in $\trip{1}{2}=e$ and get $\trip{2}{4}=e$; we also substitute $\ug{4}=\ug{3}\ug{5}\ug{3}$ in $\trip{3}{4}=e$ and get $\trip{3}{5}=e$. In a similar way, we get the relations $\trip{2}{6}=e$, $\trip{3}{6}=e$, $\trip{4}{5}=e$, and $\trip{5}{6}=e$ as well. 

\begin{figure}[ht]
\begin{center}
\scalebox{0.2}{\includegraphics{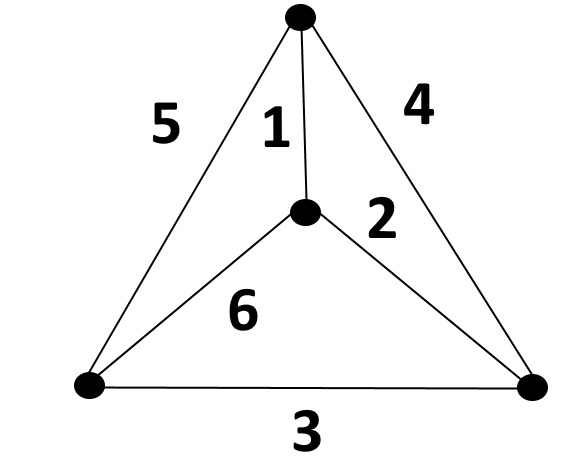}}
\end{center}
\setlength{\abovecaptionskip}{-0.15cm}
\caption{The dual graph related to $T_{(4)}$}\label{dualT4}
\end{figure}

We write all relations that we now have in group $\Ggal$ with generators $\set{\ug{i} | i=1,\dots,6}$:
	\uGammaSq{Tetra,gamma-sq-1}{1}{2}{3}{4}{5}{6}
	\begin{align}	\label{tripT4}
	&\trip{1}{2} = \trip{1}{4} = \trip{1}{5} = \trip{1}{6} = \trip{2}{3} = \trip{2}{4}= \\
        &\trip{2}{6} =\trip{3}{4} = \trip{3}{5}= \trip{3}{6} =\trip{4}{5} = \trip{5}{6} = e \nonumber
	\end{align}	
	\begin{align}	\label{commT4}
	\comm{1}{3} = \comm{2}{5} = \comm{4}{6} =e 
	\end{align}
	\begin{align}\label{cyclicT4}
	\ug{4}=\ug{1}\ug{2}\ug{1}, \ \ \ug{5}=\ug{3}\ug{4}\ug{3}, \ \ \ug{6}=\ug{2}\ug{3}\ug{2}, \ \  \ug{6}=\ug{1}\ug{5}\ug{1}.
	\end{align}
	We construct the dual graph of the edges in $T_{(4)}$, see Fig. \ref{dualT4}. 

 We can see that the relations in (\ref{tripT4}) and (\ref{commT4}) fulfill the following conditions in the graph:
 \begin{center}
\begin{eqnarray}
  [\G_i,\G_j]=e, \ \  \mbox{if $i,j$ are disjoint edges;}\label{com} \\
  \langle \G_i,\G_j \rangle=e, \ \  \mbox{otherwise,}\nonumber
\end{eqnarray}
\end{center}
and the relations in (\ref{cyclicT4}) correspond to the cycles in the graph. 
Now, by substituting $\ug{3}=\ug{4}\ug{5}\ug{4}$ in $\comm{1}{3}=e$, we get
$[\ug{1},\ug{4}\ug{5}\ug{4}] = e$. Then we substitute $\ug{5}=\ug{1}\ug{6}\ug{1}$ in $\comm{2}{5}=e$ and get
$[\ug{2},\ug{1}\ug{6}\ug{1}] = e$. In a similar way, we get the relations $[\ug{3},\ug{5}\ug{6}\ug{5}] = e$ and $[\ug{4},\ug{2}\ug{3}\ug{2}] = e$ as well. These four resulting relations correspond to the four sets of three lines each meeting at a point in the graph.
By Theorem 2.3 in \cite{RTV} we have $\Ggal \cong S_4$ (where 4 is the number of the vertices in the dual graph in Fig. \ref{dualT4}). Therefore, by (\ref{M-T}),  $\pi_1(\Xgal)$  is trivial.


\subsection{The double tetrahedron $D(T_{(4)})$}\label{DoubleTetra-sec}

In this subsection we study the fundamental group of the Galois cover of the double tetrahedron $D(T_{(4)})$. There is a nice way to construct it: we just take a union of two tetrahedrons of type $T_{(4)}$ with no bases and glue them along the edges that surround those missing bases. 
The two pieces that are being glued are depicted in Fig. \ref{twopieces} and the double tetrahedron can be seen in Fig. \ref{DoubleTetra}.

\begin{figure}[ht]
\begin{center}
\scalebox{0.50}{\includegraphics{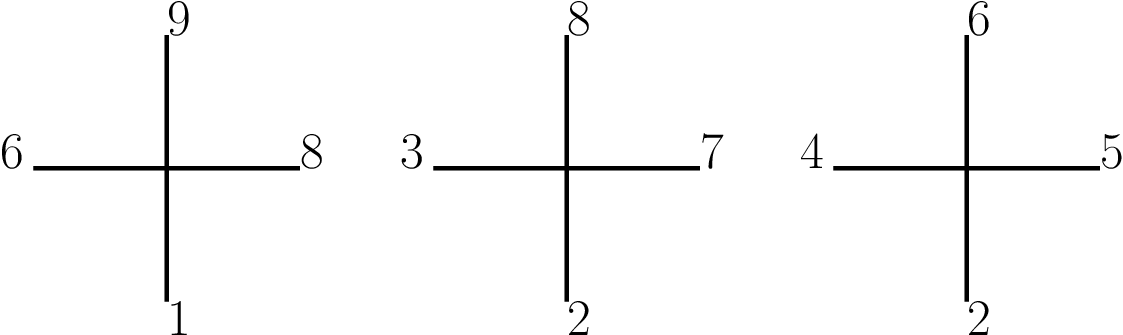}}
\end{center}
\setlength{\abovecaptionskip}{-0.15cm}
\caption{The singularities of multiplicity $4$ in $D(T_{(4)})$}\label{4-points}
\end{figure}

Vertices $V_1$ and $V_2$ are inner singularities of multiplicity $3$. The identification along edges 2, 6, and 8 creates three vertices $V_3$, $V_4$, and $V_5$, which are inner singularities of multiplicity $4$, see Fig. \ref{4-points}.

For the calculations in this section, we recall some definitions from \cite{Law} and  \cite{RTV}. Let $t \geq 0$ and $n$ be natural numbers. Let $U$ be a set consisting of $t$ elements. The group $A_{t,n}$  is generated by $n^2 t$ elements $u_{x,y}$, where $u \in U$ and $x,y \in \{1,\dots,n\}$, satisfying the following relations for any $u,u' \in U$ and any $x,y,z$:
  $u_{x,x} = e$, \ $u_{x,y}u_{y,z} = u_{x,z} = u_{y,z}u_{x,y}$, and $[u_{x,y},u_{w,z}'] = e$ (for distinct $x,y,w$, and $z$).
  
\begin{figure}[ht]
\begin{center}
\scalebox{0.80}{\includegraphics{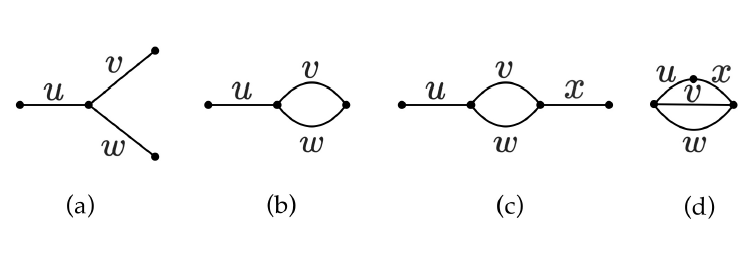}}
\end{center}
\setlength{\abovecaptionskip}{-0.15cm}
\caption{The four types of relations}\label{4cases}
\end{figure}

Given a connected graph $T$ (not necessarily simple, but without loops), the Coxeter group $C(T)$ has generators induced by the edges of $T$, and  three types of relations: $u^2 = e$ for every generator; $(uv)^2=e$ if the edges corresponding to $u,v$ are disjoint; and $(uv)^3 = e$ if $u,v$ meet in a vertex. We define the group $C_Y(T)$ as the quotient of $C(T)$ with respect to four types of relations: $[w u w, v]  =  e$ if $u,v,w$ are as in Fig. \ref{4cases}(a),   $\langle w u w,{v} \rangle  =  e$ if $u,v,w$ are as in Fig. \ref{4cases}(b),  $[w u w, v x v]  =  e$ if $x,u,v,w$ are as in Fig. \ref{4cases}(c), and $\langle w u w ,v xv \rangle  =  e$ if $x,u,v,w$ are as in Fig. \ref{4cases}(d).  

Now, the main result of \cite{RTV} is that for any planar graph $T$, simple or not, the Coxeter group $C_Y(T)$ is isomorphic to $S_n \ltimes A_{t,n}$,  where $n$ is the number of vertices in $T$ and $t$ is the rank of $\pi_1(T)$. The symmetric group acts on $A_{t,n}$, by its action on the indices. In this work we use the case $n \geq 5$ that is verified in  \cite[Theorem~5.7]{RTV}. 

The motivation for using this isomorphism is because in cases we have studied,  there was a  well-defined surjection from $C_Y(T) \rightarrow \widetilde{G}$ that sends the generator corresponding to an edge $i$ in $T$ to the generator $\Gamma_i$. 
So, what we do is take the defining relations of  $\widetilde{G}$ and pull them back to $S_n \ltimes A_{t,n}$. 

In the proof of Theorem \ref{DoubleTetGal}, we calculate the fundamental group $\pi_1(\Xgal)$ 
of the Galois cover of the double tetrahedron $D(T_{(4)})$.

\begin{theorem}\label{DoubleTetGal}
The fundamental group $\pi_1(\Xgal)$  of the Galois cover of the double tetrahedron $D(T_{(4)})$ is $\mathbb{Z}_2^4$.
\end{theorem}

\begin{verbatim}
Proof.
\end{verbatim}
There are nine edges in $D(T_{(4)})$, see Fig. \ref{DoubleTetra}. Projecting the degeneration from Fig. \ref{twopieces} onto $\mathbb{CP}^2$, we get a branch curve $S_0$, which is composed of the projections of the nine edges, and this means that $S_0$ is a line arrangement. 

There are five vertices in $D(T_{(4)})$; all of them are inner singularities with the following properties: 
Vertices $V_1$ and $V_2$ are singularities of multiplicity~$3$, and vertices $V_3$, $V_4$, and $V_5$ are singularities of multiplicity~$4$.

The branch curve $S$ has degree 18, which arises from the doubling in the regeneration. The generators in group $G$ are $\Gamma_{1},\Gamma_{1'},\dots,\Gamma_9,\Gamma_{9'}$, and we use Lemmas \ref{lemma-3} and \ref{lemma-4} for the relations coming from the singularities with multiplicity $3$ and $4$, respectively. 

We start with vertex $V_1$ and its induced relations:
\uThreePointInner{1}{3}{5}{9}{DoubleTetra,vert1}
We have the same relations to vertex $V_2$, but with indices $i=1,j=4,k=7$.

Now we give the induced relations by vertex $V_3$: 
\uFourPointInner{3}{1}{6}{8}{9}{DoubleTetra,vert3}
Similarly, we have relations that are induced by vertices $V_4$ (with indices $i={2},j={3},k={7},l={8}$) and $V_5$ (with indices $i={2},j={4},
k={5},l={6}$).

We also have the following commutative relations \eqref{DoubleTetra,parasit-1-2}-\eqref{DoubleTetra,parasit-7-9} of type (e) from Theorem \ref{vkthm} and the projective relation \eqref{DoubleTetra,proj}:
\uParasit{1}{2}{DoubleTetra,parasit-1-2}
\uParasit{1}{3}{DoubleTetra,parasit-1-3}
\uParasit{1}{5}{DoubleTetra,parasit-1-5}
\uParasit{2}{9}{DoubleTetra,parasit-2-9}
\uParasit{3}{4}{DoubleTetra,parasit-3-4}
\uParasit{3}{6}{DoubleTetra,parasit-3-6}
\uParasit{4}{8}{DoubleTetra,parasit-4-8}
\uParasit{4}{9}{DoubleTetra,parasit-4-9}
\uParasit{5}{7}{DoubleTetra,parasit-5-7}
\uParasit{5}{8}{DoubleTetra,parasit-5-8}
\uParasit{6}{7}{DoubleTetra,parasit-6-7}
\uParasit{7}{9}{DoubleTetra,parasit-7-9}
\uProjRel{DoubleTetra,proj}{9}{8}{7}{6}{5}{4}{3}{2}{1}

Because we are interested in the fundamental group of the Galois cover of $D(T_{(4)})$, we calculate now $\Ggal$. 
Vertices $V_1$ and $V_2$ are both singularities of multiplicity $3$, so we have $ \Gamma_1=\Gamma_{1'}$ and $ \Gamma_3=\Gamma_{3'}$ by Lemma \ref{3pt-in-bigmid}. Using these equalities, we can produce from the two relations (\ref{4pt-rels-9}) and (\ref{4pt-rels-10}), which are induced by $V_4$ (indices $i=2,j=3,k=7,l=8$), the following two relations:
\begin{eqnarray}
\Gamma_{2'}\Gamma_3\Gamma_{2'} = \Gamma_8\Gamma_{7'}\Gamma_8 \label{num1}\\
\Gamma_{2'}\Gamma_3\Gamma_{2'} = \Gamma_8\Gamma_{7'}\Gamma_7{\Gamma_{7'}}\Gamma_8. \label{num2}
\end{eqnarray}
We equate \eqref{num1} and \eqref{num2}
$$\Gamma_8\Gamma_{7'}\Gamma_8 = \Gamma_8\Gamma_{7'}\Gamma_7{\Gamma_{7'}}\Gamma_8$$
and get $$\Gamma_{7'} = \Gamma_7.$$
This result is applied in Lemma \ref{3pt-in-bigmid} for vertex $V_2$, giving us $\Gamma_4=\Gamma_{4'}$. 
Now we take again the two relations  (\ref{4pt-rels-9}) and (\ref{4pt-rels-10}) that correspond to vertex $V_5$ (indices $i=2,j=4,k=5,l=6$) and substitute $\Gamma_4=\Gamma_{4'}$ to get:
$$\Gamma_{2'}\Gamma_4\Gamma_{2'} = \Gamma_6\Gamma_{5'}\Gamma_6$$
$$\Gamma_{2'}\Gamma_4\Gamma_{2'} = \Gamma_6\Gamma_{5'}\Gamma_5{\Gamma_{5'}}\Gamma_6.$$
We equate them and get $\Gamma_5=\Gamma_{5'}$ as well.
Then, applying Lemma \ref{3pt-in-bigmid} on vertex $V_1$, with the equality $\Gamma_5=\Gamma_{5'}$, we get $\Gamma_9=\Gamma_{9'}$.

Next, let us look at $\Gamma_{2'}\Gamma_3\Gamma_{2'} = \Gamma_8\Gamma_{7'}\Gamma_8$ from \eqref{num1}. Because we have $\Gamma_{7'} = \Gamma_7$, the relation becomes $\Gamma_{2'}\Gamma_3\Gamma_{2'} = \Gamma_8\Gamma_{7}\Gamma_8$. Then we use  $\langle \Gamma_{2'}, \Gamma_3  \rangle=e$ (coming from (\ref{4pt-rels-1}) for $i=2,j=3$), to rewrite the relation as $\Gamma_{3}\Gamma_{2'}\Gamma_{3} = \Gamma_8\Gamma_{7}\Gamma_8$, and therefore we get: 
\begin{equation}\label{2a}
\Gamma_{2'} = \Gamma_3 \Gamma_8\Gamma_7\Gamma_8\Gamma_3.
\end{equation}
Similarly, we get the equality:
\begin{equation}\label{2b}
\Gamma_{2'} = \Gamma_4\Gamma_6\Gamma_5\Gamma_6\Gamma_4.
\end{equation}
From (\ref{4pt-rels-11}) (with indices $i=2,j=3,k=7,l=8$ and also $i=2,j=4,k=5,l=6$), we can get the following two equalities:
\begin{equation}\label{8eliminate}
\Gamma_{8'} = \Gamma_8\Gamma_7\Gamma_2\Gamma_3\Gamma_2\Gamma_7\Gamma_8
\end{equation}
\begin{equation}\label{6eliminate}
\Gamma_{6'} = \Gamma_6\Gamma_5\Gamma_2\Gamma_4\Gamma_2\Gamma_5\Gamma_6.
\end{equation}

Our next step is to use $\Gamma_i=\Gamma_{i'}$ ($i=1,3,4,5,7,9$), and also expressions \eqref{2a}-\eqref{6eliminate}, to eliminate the generators $\Gamma_{2'}$, $\Gamma_{6'}$, and $\Gamma_{8'}$.  This elimination is accompanied by the placement of those expressions in all relations in $\Ggal$, and then many group reductions.  New relations appear, for example  $\trip{1}{7} =e$,  $[\Gamma_3, \Gamma_7]=e$, and $[\ug{1},\ug{8}\ug{7}\ug{8}]=e$. We explain here how to derive them; the other new relations that did not exist before can be produced in a similar way as well.  
First,  we show how to get $\trip{1}{7} =e$: relation (\ref{3pt-rels-3}) (with $i=1,j=4,k=7$) can be rewritten as $\Gamma_{4} = \Gamma_1\Gamma_7\Gamma_1$, which can be substituted into $\trip{1}{4}=e$ (from (\ref{3pt-rels-2}), again for $i=1,j=4,k=7$), with the resulting relation being $\trip{1}{7} =e$. Second, we show how to get   $[\Gamma_3, \Gamma_7]=e$. Starting with $\trip{2'}{3}=e$ ((\ref{4pt-rels-1}) for $i=2,j=3$), we substitute \eqref{2a} inside, and get $\langle \Gamma_3, \Gamma_8\Gamma_7\Gamma_8  \rangle=e$. Then we substitute \eqref{2a} inside $\comm{2'}{8}=e$ ((\ref{4pt-rels-3}) for $i=2,l=8$), and get $[\Gamma_3 \Gamma_8\Gamma_7\Gamma_8\Gamma_3, \Gamma_8]=e$. We use the relation $\langle \Gamma_3, \Gamma_8\Gamma_7\Gamma_8  \rangle=e$ to rewrite it as $[\Gamma_8\Gamma_7\Gamma_8\Gamma_3 \Gamma_8\Gamma_7\Gamma_8, \Gamma_8]=e$. Because we have $\trip{7}{8} = e$ ((\ref{4pt-rels-2}) for $k=7,l=8$), we get $[\Gamma_3, \Gamma_7]=e$.
Third, we substitute \eqref{2a} in  $\comm{1}{2'}=e$ (coming from \eqref{DoubleTetra,parasit-1-2}) and get $[\Gamma_1,\Gamma_3\Gamma_8\Gamma_7\Gamma_8\Gamma_3]=e$, which can be simplified to $[\Gamma_1,\Gamma_8\Gamma_7\Gamma_8]=e$ by \eqref{DoubleTetra,parasit-1-3}. 

At this stage, after eliminating  $\Gamma_{2'}$, $\Gamma_{6'}$, and $\Gamma_{8'}$, we write the presentation for group $\Ggal$ with the generators $\Gamma_i$, \ $i=1, \dots, 9$ and the following relations:
\uGammaSq{DoubleTetra,gamma-sq1}{1}{2}{3}{4}{5}{6}{7}{8}{9}
\begin{align}
	\trip{1}{6} = \trip{1}{7} = \trip{2}{3} = \trip{2}{7} = \trip{3}{5} = \trip{3}{8} &= \label{simp1}\\
	\trip{3}{9} =\trip{4}{6} = \trip{4}{7}= \trip{5}{6} = \trip{6}{9} = \trip{7}{8} &= e \nonumber
\end{align}
\begin{align}
    \comm{1}{3} = \comm{2}{6} = \comm{3}{4} = \comm{3}{6} = \comm{3}{7} = \comm{5}{7}  &= \label{simp2}\\ 
    \comm{6}{7} =  \comm{6}{8} = \comm{7}{9} &= e \nonumber
      \end{align}
\begin{align}
    [\ug{1},\ug{8}\ug{7}\ug{8}]=[\ug{2},\ug{4}\ug{7}\ug{4}]=[\ug{2},\ug{5}\ug{3}\ug{5}]&=\label{simp3}\\
      [\ug{6},\ug{5}\ug{9}\ug{5}]=[\ug{6},\ug{1}\ug{4}\ug{1}]=[\ug{8},\ug{3}\ug{9}\ug{3}] &= e \nonumber
\end{align}
\begin{align} \label{simp4}
\ug{7}=\ug{1}\ug{4}\ug{1}, \ \ \ug{3}=\ug{5}\ug{9}\ug{5}, \ \ \ug{3}\ug{8}\ug{7}\ug{8}\ug{3}
=\ug{4}\ug{6}\ug{5}\ug{6}\ug{4},
\end{align}
and
\begin{align} \label{left1}
	\trip{1}{4} = \trip{1}{8} = \trip{2}{4} = \trip{2}{5} = \trip{5}{9} = \trip{8}{9} = e
\end{align}
\begin{align}
    \comm{1}{2} = \comm{1}{5} = \comm{1}{9} = \comm{2}{8} = \comm{2}{9} &= \label{left2}\\
    \comm{4}{5} = \comm{4}{8} =  \comm{4}{9} = \comm{5}{8} &= e.\nonumber
      \end{align}
We still have the simplified projective relation,
\begin{equation}\label{left3}
\Gamma_{8'}\Gamma_8\Gamma_{6'}\Gamma_6\Gamma_{2'}\Gamma_2 = e,
\end{equation}
and we will deal with it further in the following paragraphs.
%
%
%

Now we eliminate the generators $\Gamma_7, \Gamma_3$, and $\Gamma_6$, using the equalities from \eqref{simp4}.
We start with the eliminations of $\ug{7}=\ug{1}\ug{4}\ug{1}$ and $\ug{3}=\ug{5}\ug{9}\ug{5}$, which turn to be  quite immediate.
After this, we substitute $\ug{7}=\ug{1}\ug{4}\ug{1}$ and  $\ug{3}=\ug{5}\ug{9}\ug{5}$ in
$\ug{3}\ug{8}\ug{7}\ug{8}\ug{3}=\ug{4}\ug{6}\ug{5}\ug{6}\ug{4}$ (from \eqref{simp4}) to get
$\ug{6}=\ug{9}\ug{8}\ug{1}\ug{8}\ug{9}$. Then we also eliminate $\Gamma_6$.
The relations that appear in \eqref{simp1}-\eqref{simp3} become redundant. The following example presents a relation that is redundant (the remaining relations are similarly redundant and also include the same process, which is relatively immediate): take $\trip{1}{6} =e$ and substitute inside the expression $\ug{6}=\ug{9}\ug{8}\ug{1}\ug{8}\ug{9}$ to get $\langle\ug{1},\ug{9}\ug{8}\ug{1}\ug{8}\ug{9}\rangle=e$. By $ \comm{1}{9} =e$, we get $\langle\ug{1},\ug{8}\ug{1}\ug{8}\rangle=e$, which is equivalent to $\trip{1}{8}=e$. 

\begin{figure}[H]
\begin{center}
\scalebox{0.30}{\includegraphics{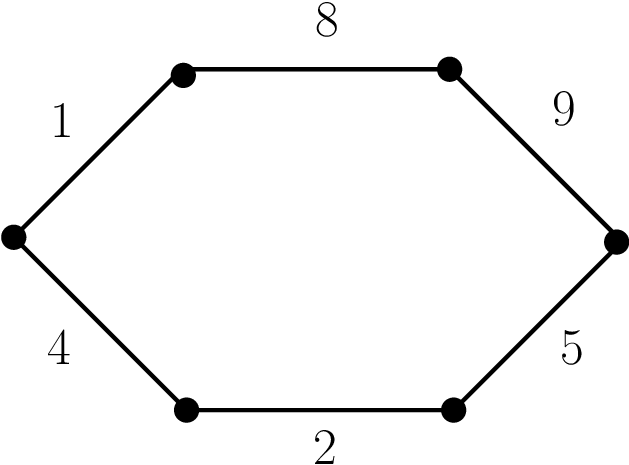}}
\end{center}
\setlength{\abovecaptionskip}{-0.15cm}
\caption{Graph $ T $ related to $ C_Y(T) $}\label{meshushe}
\end{figure}

We get a much more simplified presentation of $\Ggal$, with generators
$\ug{1}, \ug{2}, \ug{4}, \ug{5}, \ug{8}, \ug{9}$, relations \eqref{left1} and \eqref{left2}, the relation
\uGammaSq{DoubleTetra,gamma-sq}{1}{2}{4}{5}{8}{9}
and the relation we will obtain from \eqref{left3} after performing the substitutions mentioned above (let us denote it as $Proj$).
In our case, without relation \eqref{left3}, this presentation is exactly the Coxeter group $ C_Y(T) $, where $ T $ is the graph depicted in Fig. \ref{meshushe}.
When $ t=1 $ (i.e., there is one cycle in $T$), the group $ A_{1,n} $ is the Abelian group generated by $\{ u_{i,j} \; \mid \; 1\le i,j \le n \} $ with the relations $ u_{i,i}=e $, $u_{i,j}=u_{j,i}^{-1}$,
$u_{i,k}u_{k,j}=u_{i,j}$; it is obviously isomorphic to $ \mathbb{Z}^{n-1} $ (see \cite[Example 5.3]{RTV}).
Here we have  $t=1$ and $n=6$, and therefore  
$$C_Y(T) \cong S_6 \ltimes A_{1,6} = S_6 \ltimes \mathbb{Z}^5.$$
We can see that
$$\Ggal = \frac{C_Y(T)}{H},$$ where $H$ is the normal subgroup generated by $ Proj $ in $C_Y(T)$.

\begin{figure}[ht]
\begin{center}
\scalebox{0.25}{\includegraphics{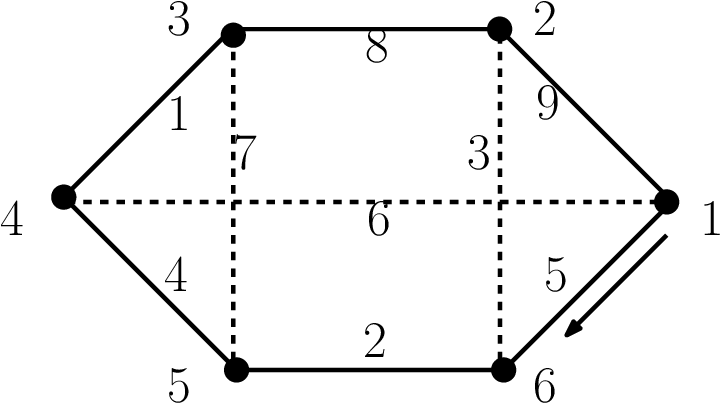}}
\end{center}
\setlength{\abovecaptionskip}{-0.15cm}
\caption{The graph with assignments of elements}\label{meshushe-CY}
\end{figure}

Now we assign to each element $\Gamma_i$ ($i=1,2,4,5,8,9$) in $C_Y(T)$  an element in $S_6 \ltimes \mathbb{Z}^5$, see Fig. \ref{meshushe-CY}. We choose the spanning tree $ T_0 $ in the proof of \cite[Theorem 6.1]{RTV} to be all the edges but $ \Gamma_5 $:
$$\Gamma_1=(3 \; 4), \  \Gamma_2=(5 \; 6), \ \Gamma_4=(4 \; 5), \ \Gamma_5=(1 \; 6)u_{1,6},
\ \Gamma_8=(2 \; 3), \ \Gamma_9=(1 \; 2).$$

Because the action of $S_n$ on $A_{t,n}$ is defined by $\sigma u_{i,j}=u_{\sigma(i),\sigma(j)}$, we can determine the assigned elements of $\Gamma_3, \Gamma_6, \Gamma_7,
\Gamma_{2'}, \Gamma_{6'}$, and $\Gamma_{8'}$ as follows:

\begin{description}
  \item[$\ug{3}$:] \ \ $\ug{3}=\ug{9}\ug{5}\ug{9}=(1 \; 2)(1 \; 6)u_{1,6}(1 \; 2)=(2 \; 6)u_{2,6}$
  \item[$\ug{6}$:] \ \ $\ug{6}=\ug{9}\ug{8}\ug{1}\ug{8}\ug{9}=(1 \; 2)(2 \; 3)(3 \; 4)(2 \; 3)(1 \; 2)=(1 \; 4)$
  \item[$\ug{7}$:] \ \ $\ug{7}=\ug{1}\ug{4}\ug{1}=(3 \; 4)(4 \; 5)(3 \; 4)=(3 \; 5)$
  \item[$\Gamma_{8'}$:] \ \
$\Gamma_2\Gamma_3\Gamma_2=(5 \; 6)(2 \; 6)u_{2,6}(5 \; 6)=(2 \; 5)u_{2,5} \Longrightarrow$\\
$\Gamma_7\Gamma_2\Gamma_3\Gamma_2\Gamma_7=(3 \; 5)(2 \; 5)u_{2,5} (3 \; 5)=(2 \; 3)u_{2,3} \Longrightarrow$\\
$\Gamma_{8'} = \Gamma_8\Gamma_7\Gamma_2\Gamma_3\Gamma_2\Gamma_7\Gamma_8=(2 \; 3)(2 \; 3)u_{2,3}(2 \; 3)
=(2 \; 3)u_{3,2}=(2 \; 3)u_{2,3}^{-1}$
  \item[$\Gamma_{6'}$:] \ \
$\Gamma_2\Gamma_4\Gamma_2=(5 \; 6)(4 \; 5)(5 \; 6)=(4 \; 6) \Longrightarrow$\\
$\Gamma_5\Gamma_2\Gamma_4\Gamma_2\Gamma_5=(1 \; 6)u_{1,6}(4 \; 6)(1 \; 6)u_{1,6}=
(1 \; 6)u_{1,6}(4 \; 1 \; 6)u_{1,6}=(1 \; 6)(4 \; 1 \; 6)u_{6,4}u_{1,6}=(1 \; 4)u_{1,4} \Longrightarrow$\\
$\Gamma_{6'} = \Gamma_6\Gamma_5\Gamma_2\Gamma_4\Gamma_2\Gamma_5\Gamma_6=(1 \; 4)(1 \; 4)u_{1,4}(1 \; 4)=
(1 \; 4)u_{4,1}=(1 \; 4)u_{1,4}^{-1}$
   \item[$\Gamma_{2'}$:] \ \
$\Gamma_8\Gamma_7\Gamma_8=(2 \; 3)(3 \; 5)(2 \; 3)=(2 \; 5) \Longrightarrow$\\
$\Gamma_{2'} = \Gamma_3\Gamma_8\Gamma_7\Gamma_8\Gamma_3=(2 \; 6)u_{2,6}(2 \; 5)(2 \; 6)u_{2,6}=
(2 \; 6)(2 \; 5)(2 \; 6)u_{5,2}u_{2,6}=(5 \; 6)u_{5,6}$.
\end{description}

Next, we calculate the value of $ Proj $ in $S_6 \ltimes \mathbb{Z}^5 $:
\begin{equation}\label{proj-new}
\begin{split}
e=\Gamma_{8'}\Gamma_8\Gamma_{6'}\Gamma_6\Gamma_{2'}\Gamma_2=
(2 \; 3)u_{2,3}^{-1}(2 \; 3)(1 \; 4)u_{1,4}^{-1}(1 \; 4)(5 \; 6)u_{5,6}(5 \; 6)=
u_{3,2}^{-1}u_{4,1}^{-1}u_{6,5}=\\
u_{2,3}u_{1,4}u_{5,6}^{-1}=u_{1,4}u_{2,3}u_{5,6}^{-1}=
(u_{1,2}u_{2,3}u_{3,4})u_{2,3}u_{5,6}^{-1}=u_{1,2}u_{2,3}^2u_{3,4}u_{5,6}^{-1}.
\end{split}
\end{equation}
Each element of the form $u_{i,i+1}$  in $S_6 \ltimes \mathbb{Z}^5$ fulfills $u^2_{i,i+1} \in H$ because we have:
$$u_{5,6}^2 = u_{5,6}u_{5,6} =  (u_{1,2}u_{2,3}^2u_{3,4}u_{5,6})u_{1,2}^{-1}u_{2,3}^{-2}u_{3,4}^{-1}u_{5,6} = ((5 \; 6)Proj(5 \; 6)){Proj}^{-1} \in H $$
$$\Longrightarrow$$
$$u_{4,5}^2=(1 \; 2\; 3\; 4\; 5\; 6)u_{5,6}^2(6 \; 5\; 4\; 3\; 2\; 1) \in H $$
$$\Longrightarrow$$
$$u_{5-k, 6-k}^2 =(1 \; 2\; 3\; 4\; 5\; 6)^k u_{5,6}^2 (6 \; 5\; 4\; 3\; 2\; 1)^k \in H.$$
Thus, 
$$\Ggal = \frac{S_6 \ltimes \mathbb{Z}^5}{H}=\frac{S_6 \ltimes \mathbb{Z}_2^5}{N},$$
where $N$ is the normal subgroup of $S_6 \ltimes \mathbb{Z}_2^5$ generated by the image of $Proj$ in this group, i.e., by $u_{1,2}u_{3,4}u_{5,6} \in S_6 \ltimes \mathbb{Z}_2^5$.

We now show that the group generated by $u_{1,2}u_{3,4}u_{5,6}$ in $S_6 \ltimes \mathbb{Z}_2^5$ is normal (and thus equal to $N$).
It is enough to show that $\sigma^{-1} u_{1,2}u_{3,4}u_{5,6} \sigma \in \langle
u_{1,2}u_{3,4}u_{5,6} \rangle$, for the standard generators $\sigma=(i\; i+1)$ of $S_6$, \ $1 \leq i \leq 5$.

\begin{description}
  \item[Conjugation 1:] \
$(1 \; 2)Proj(1 \; 2)=(1 \; 2)u_{1,2}u_{3,4}u_{5,6}(1 \; 2)=u_{2,1}u_{3,4}u_{5,6}=
u_{1,2}u_{3,4}u_{5,6}=Proj$
  \item[Conjugation 2:] \
$(2 \; 3)Proj(2 \; 3)=(2 \; 3)u_{1,2}u_{3,4}u_{5,6}(2 \; 3)=u_{1,3}u_{2,4}u_{5,6}=
(u_{1,2}u_{2,3})(u_{2,3}u_{3,4})u_{5,6}=u_{1,2}u_{3,4}u_{5,6}=Proj$
  \item[Conjugation 3:] \
$(3 \; 4)Proj(3 \; 4)=(3 \; 4)u_{1,2}u_{3,4}u_{5,6}(3 \; 4)=u_{1,2}u_{4,3}u_{5,6}=
u_{1,2}u_{3,4}u_{5,6}=Proj$
  \item[Conjugation 4:] \
$(4 \; 5)Proj(4 \; 5)=(4 \; 5)u_{1,2}u_{3,4}u_{5,6}(4 \; 5)=u_{1,2}u_{3,5}u_{4,6}=
u_{1,2}(u_{3,4}u_{4,5})(u_{4,5}u_{5,6})=u_{1,2}u_{3,4}u_{5,6}=Proj$
  \item[Conjugation 5:] \
$(5 \; 6)Proj(5 \; 6)=(5 \; 6)u_{1,2}u_{3,4}u_{5,6}(5 \; 6)=u_{1,2}u_{3,4}u_{6,5}=
u_{1,2}u_{3,4}u_{5,6}=Proj$.
\end{description}

Therefore $N$ is generated by $u_{1,2}u_{3,4}u_{5,6} \in S_6 \ltimes \mathbb{Z}_2^5$, so by linear algebra 
$$\Ggal=S_6 \ltimes \mathbb{Z}_2^4,$$ 
and by Equation (\ref{M-T}): $$\pi_1(\Xgal)=\mathbb{Z}_2^4.$$



\section{Conclusion}\label{conc}
In Subsection \ref{Tetra-sec}, we determined that $\pi_1(\Xgal)$  is trivial, which means that the Galois cover of the tetrahedron is simply-connected. In Subsection \ref{DoubleTetra-sec}, we  determined that $\mathbb{Z}_2^4$ is the corresponding group to the double tetrahedron.

\begin{corollary}
The surfaces $T_{(4)}$ and $D(T_{(4)})$ are located in different connected components in the moduli space of surfaces.
\end{corollary}

We are now motivated to continue the inductive investigation of surfaces having non-planar degeneration. We will probably have complicated algebraic computations in groups as well as daunting complications while resolving  high-multiplicity singular points (such as the singularities of multiplicity $k$ presented in this work).  Lemma \ref{3pt-in-bigmid} gave us a rule about singularities of multiplicity $3$ that eases group calculations, compared to singularities of multiplicity $4$ that appear in gluing two pieces, for example, as happens for  $ D(T_{(4)}) $ or for any two Zappatic degenerations (see e.g., \cite{RkRk}), and for which we still have no similar lemma. We could determine the group very precisely in both cases and get the above corollary. 
We will need to develop strategies for attacking various singularities of multiplicity $k$ and try to determine the general forms of groups.

\bmhead{Acknowledgments}

 I thank the anonymous referees for their mathematical suggestions, and also for the suggestion for a much more clear and orderly structure of the paper, which was implemented. 

 Thanks to Uriel Sinichkin (Tel Aviv University) and Cheng Gong (Soochow University) for efficient discussions.



%

\end{document}